\documentclass{amsart}
\usepackage{amssymb}
\usepackage{amscd}

\newcommand\Zs{{\mathbb Z}}
\newcommand\Ns{{\mathbb N}_0}
\newcommand\ks{{\Bbbk}}
\newcommand\End{\mathop {\fam 0 End}\nolimits}

\newcommand\Der{\mathop {\fam 0 Der}\nolimits}

\newcommand\cha{\mathop {\fam 0 char}\nolimits}

\begin{document}

\title[SUBALGEBRAS IN POSITIVE CHARACTERISTIC]
{SUBALGEBRAS OF THE POLYNOMIAL ALGEBRA IN POSITIVE CHARACTERISTIC
AND THE JACOBIAN}
\author{Alexey~V.~Gavrilov}
\email{gavrilov19@gmail.com}
\maketitle
{\small
Let $\ks$ be a field of characteristic $p>0$ and
$R$ be a subalgebra of $\ks[X]=\ks[x_1,\dots,x_n]$. Let
$J(R)$ be the ideal in $\ks[X]$ defined by
$J(R)\Omega_{\ks[X]/\ks}^n=\ks[X]\Omega_{R/\ks}^n$. It is shown that if
it is a principal ideal then $J(R)^q\subset R[x_1^p,\dots,x_n^p]$,
where $q=\frac{p^n(p-1)}{2}$.}

\section{INTRODUCTION}

Let $\ks$ be a field and $\ks[X]=\ks[x_1,\dots,x_n]$ be the polynomial ring.
For the most of the paper the number $n\ge 1$ and the basis
$X=\{x_1,\dots,x_n\}$ are fixed.
Let  $R$ be a subalgebra of $\ks[X]$.
Denote by $J(R)$ the ideal in $\ks[X]$ generated by the Jacobians
of sets of elements of $R$ (a formal definition will be given below).
The main result is the following theorem.

{\bf Theorem}
{\it Let $\ks$ be a field of characteristic $p>0$ and
$R$ be a subalgebra of $\ks[X]$. If $J(R)$ is a principal ideal then
$$J(R)^{q}\subset R[X^p],$$
where $R[X^p]=R[x_1^p,\dots,x_n^p]$ and $q=\frac{p^n(p-1)}{2}$.
}

Presumably the statement holds for a nonprincipal ideal as well.
However, a proof of this conjecture probably requires another technique.

{\bf Corollary }
{\it
Let $\ks$ be a field of characteristic $p>0$ and
$R$ be a subalgebra of $\ks[X]$. If $J(R)=\ks[X]$ then
$$R[X^p]=\ks[X].$$
}

Nousiainen [2] proved this in the case $R=\ks[F]=\ks[f_1,\dots,f_n]$
(see also [1]). In this case the condition $J(R)=\ks[X]$ is
equivalent to $j(F)\in\ks^{\times}$, where $j(F)$ is the Jacobian of
the polynomials $f_1,\dots,f_n$. Thus, Nousiainen's result is a
positive characteristic analogue of the famous Jacobian conjecture:
if $\cha(\ks)=0$ and $j(F)\in\ks^{\times}$ then $\ks[F]=\ks[X]$. The
zero characteristic analogue of Corollary is obviously false.
(Consider, for example, the subalgebra $R=\ks[t-t^2,t-t^3]$ of
$\ks[t]$. Then $R\neq\ks[t]$ but $J(R)=\ks[t]$).

Nousiainen's method is based on properties of the derivations
$\frac{\partial}{\partial f_i}\in\Der(\ks[X])$, which are the
natural derivations of $\ks[F]$ extended to $\ks[X]$. Probably his
method can be applied to a more general case as well. However, our
approach is based on the calculation of generalized Wronskians of a
special kind. This calculation may be an interesting result in
itself.

\section{PRELIMINARIES AND NOTATION}

An element $\alpha=(\alpha_1,\dots,\alpha_n)$ of $\Ns^n$, where
$\Ns$ denotes the set of non-negative integers, is called a multiindex.
For $F\in\ks[X]^n$ and two multiindices
$\alpha,\beta\in\Ns^n$ we will use the following notation

$$|\alpha|=\sum_{i=1}^n\alpha_i,\,\alpha!=\prod_{i=1}^n\alpha_i!,\,
{\alpha\choose\beta}=\prod_{i=1}^n{\alpha_i\choose\beta_i},$$

$$X^{\alpha}=\prod_{i=1}^n x_i^{\alpha_i},\,
F^{\alpha}=\prod_{i=1}^n f_i^{\alpha_i},\,
\partial^{\alpha}=\prod_{i=1}^n \partial_i^{\alpha_i},$$
where $\partial_i=\frac{\partial}{\partial x_i}\in\Der(\ks[X]).$

In several places multinomial coefficients will appear, denoted by
${\alpha\choose\beta^1\dots\beta^k}$, where
$\alpha,\beta^1,\dots,\beta^k\in\Ns^n$ and
$\alpha=\beta^1+\dots+\beta^k$. They are defined exactly by the same
way as the binomial ones. The set of multiindices possess the
partial order; by definition, $\alpha\le\beta$ iff
$\alpha_i\le\beta_i$ for all $1\le i\le n$. For the sake of
convenience we introduce the "diagonal" intervals
$[k,m]=\{\alpha\in\Ns^n:k\le\alpha_i\le m, 1\le i\le n\}$, where
$k,m\in\Ns$.

Let $R$ be a subalgebra of $\ks[X]$. Since $\Omega_{\ks[X]/\ks}^n$
is a free cyclic module, we can make the following definition.

{\bf Definition 1}
{\it
Let $R$ be a subalgebra of $\ks[X]$,  where $\ks$ is a field.
The Jacobian ideal of $R$ is the ideal
$J(R)$ in $\ks[X]$ defined by the equality
$$J(R)\Omega_{\ks[X]/\ks}^n=\ks[X]\Omega_{R/\ks}^n,$$
where $\Omega_{R/\ks}^n$ is considered a submodule of $\Omega_{\ks[X]/\ks}^n$
over $R$.
}

The exact meaning of the words "is considered a submodule" is that we write
$\ks[X]\Omega_{R/\ks}^n$ instead of
$\ks[X]{\rm Im}(\Omega_{R/\ks}^n\to\Omega_{\ks[X]/\ks}^n)$. This is a slight
abuse of notation, because the natural $R$ - module homomorphism
$\Omega_{R/\ks}^n\to\Omega_{\ks[X]/\ks}^n$ is not injective in general.

There is also a more explicit definition. For $F\in\ks[X]^n$, the
Jacobian matrix and the Jacobian are defined by
$$JF=\left\|\frac{\partial f_i}{\partial x_j}\right\|_{1\le i,j\le n}\in
M(n,\ks[X]),\,
j(F)=\det JF\in\ks[X].$$

The module $\ks[X]\Omega_{R/\ks}^n$ is generated by
$df_1\wedge\dots\wedge df_n=j(F)dx_1\wedge\dots\wedge dx_n$, where
$F=(f_1,\dots,f_n)\in R^n$. Thus
$$J(R)=\langle\{j(F):F\in R^n\}\rangle,$$
where $\langle S\rangle$ denotes the ideal in $\ks[X]$ generated by a set $S$.
It is an easy consequence of the chain rule that the Jacobian ideal
of a subalgebra generated by $n$ polynomials is a principal ideal:
$$J(\ks[F])=\ks[X]j(F),\,F\in\ks[X]^n.$$

Clearly, the ring $\ks[X]$ is a free module over $\ks[X^p]$ of rank
$p^n$: the set of monomials $\{X^{\alpha}:\alpha\in[0,p-1]\}$ is a
natural basis of this module. This construction became important
when $\cha(\ks)=p$ (note that in this case $\ks[X^p]$ does not
depend on the choice of generators of $\ks[X]$ i.e. it is an
invariant).

{\bf Definition 2}
{\it
Let $\ks$ be a field of characteristic $p>0$ and $F\in\ks[X]^n$.
The matrix $U(F)\in M(p^n,\ks[X^p])$ is defined by
$$F^{\alpha}=\sum_{\beta\in[0,p-1]}U(F)_{\alpha\beta}X^{\beta},
\,\alpha\in[0,p-1].$$
}

\section{GENERALIZED WRONSKIANS}

In this section we compute generalized Wronskians of a special form.
The key tool for this computation is the following simple lemma.

{\bf Lemma 1}
{\it Let $R$ be a ring and $f\in R$. If $D_1,\dots,D_l\in\Der(R)$
and if $m\ge l\ge 0$ then
$$\sum_{k=0}^m{m\choose k}(-f)^{m-k}D_1\dots D_l f^k
=\begin{cases}
0 & {\rm if}\, m>l \\
l!\prod_{k=1}^lD_kf & {\rm if}\, m=l
\end{cases}
\eqno{(1)}$$
}

In the case $l=0$ there are no derivations and the formula becomes
$$\sum_{k=0}^m{m\choose k}(-f)^{m-k} f^k
=\begin{cases}
0 & {\rm if}\, m>0 \\
1 & {\rm if}\, m=0
\end{cases}
$$
which is obviously true.
\newline
{\it Proof}

Denote
$$S_{m,l}=\sum_{k=0}^m{m\choose k}(-f)^{m-k}D_l\dots D_1 f^k.$$
This sum coincides with the left hand side of (1), except for the
reverse order of the derivations.
We have $S_{0,0}=1$ and $S_{m,0}=0,\,m>0$.
The following equality can easily be verified
$$S_{m+1,l+1}=D_{l+1}S_{m+1,l}+(m+1)(D_{l+1}f)S_{m,l}.$$
By induction on $l$, for $m>l$ we have $S_{m,l}=0$, and
$$S_{l,l}=l(D_{l}f)S_{l-1,l-1}=
l!\prod_{k=1}^lD_kf.$$
\qed

Let $R$ be a ring. Denote by $R[Z_{ij}]$ the polynomial ring in the
$n^2$ indeterminates $Z_{ij},\,1\le i,j\le n$. If $h\in R[Z_{ij}]$ and
$A\in M(n, R)$ then $h(A)$ denotes the result of the substitution
$Z_{ij}\mapsto A_{ij}$.

{\bf Proposition 1} {\it For any $r\ge 1$ there exists the
homogeneous polynomial $H_r\in \Zs[Z_{ij}]$ of degree
$\frac{nr^n(r-1)}{2}$ with the following property. Let $R$ be a ring
with derivations $D_1,\dots,D_n\in\Der(R)$. If $[D_i,D_j]=0$, for
all $i,j$, then for any $F=(f_1,\dots,f_n)\in R^n$ the following
equality holds
$$\det W=H_r(JF),\eqno{(2)}$$
where $W\in M(r^n,R)$ and $JF\in M(n,R)$ are defined by
$$W_{\alpha\beta}=D^{\alpha}F^{\beta},\,\alpha,\beta\in[0,r-1];\,
(JF)_{ij}=D_if_j,\,1\le i,j\le n.$$
}

Here $JF$ is an obvious generalization of the Jacobian matrix. The
determinant $\det W$ is a generalized Wronskian of the polynomials
$F^{\beta}$.
\newline
{\it Proof}

By the Leibniz formula,
$$W_{\alpha\beta}=\sum
{\alpha\choose\theta^1\dots\theta^n}\prod_{i=1}^n(D^{\theta^i}
f_i^{\beta_i}),$$
where $\theta^1,\dots,\theta^n\in\Ns^n$ and the sum is taken over the
multiindices satisfying the equality
$\sum_{i=1}^n\theta^i=\alpha$.

Let $T\in M(r^n,R)$ be determined by
$$T_{\alpha\beta}={\beta\choose\alpha}(-F)^{\beta-\alpha},\,
\alpha,\beta\in[0,r-1].$$ Let $W^{\prime}=WT\in M(r^n,R)$. Then
$$W_{\alpha\beta}^{\prime}=\sum_{\gamma}W_{\alpha\gamma}T_{\gamma\beta}=
\sum{\alpha\choose\theta^1\dots\theta^n}\prod_{i=1}^n
S_i(\beta_i,\theta^i),\eqno{(3)}$$
where
$$S_i(m,\theta)=\sum_{k=0}^m{m\choose k}(-f_i)^{m-k}D^{\theta}f_i^k.$$

By Lemma 1, if $m\ge|\theta|$ then

$$S_i(m,\theta)
=\begin{cases}
0 & {\rm if}\, m>|\theta| \\
m!\prod_{j=1}^n(D_j f_i)^{\theta_j} & {\rm if}\, m=|\theta|
\end{cases}
$$

Thus if the product in the right hand side of (3) is not zero
then $\beta_i\le|\theta^i|$ for all $1\le i\le n$. The latter implies
$|\beta|\le|\alpha|$. It follows that if $|\alpha|<|\beta|$ then
$W_{\alpha\beta}^{\prime}=0$. In the case $|\alpha|=|\beta|$ the product
is zero unless $\beta_i=|\theta^i|,\,1\le i\le n$. Thus, if
$|\alpha|=|\beta|$ then
$$W_{\alpha\beta}^{\prime}=\beta!
\sum_{|\theta^1|=\beta_1}\dots\sum_{|\theta^n|=\beta_n}
{\alpha\choose\theta^1\dots\theta^n}\prod_{i=1}^n\prod_{j=1}^n
(D_jf_i)^{\theta^i_j}.\eqno{(4)}$$

Put the multiindices in a total order compatible with the partial
order (e.g. in the lexicographic order). Then the matrix $T$ becomes
upper triangular with the unit diagonal, hence $\det T=1$. The
matrix $W^{\prime}$ becomes block lower triangular, hence $\det
W^{\prime}$ is equal to the product of the $nr-n+1$ determinants of
the blocks. Each block determinant
$\det\left\|W_{\alpha\beta}^{\prime}\right\|_{|\alpha|=|\beta|=l}$
(where $0\le l\le nr-n$) is by (4) a homogeneous polynomial in the
variables $D_if_j$ of degree $ls_l$, where $s_l$ is the size of the
block. Clearly $s_l=\#\{\alpha\in\Ns^n:\alpha\in[0,r-1],
|\alpha|=l\}$. The determinant $\det W=\det W^{\prime}$ is then a
homogeneous polynomial of degree
$$\sum_{l=0}^{nr-n}ls_l=\sum_{\alpha\in[0,r-1]}|\alpha|=\frac{nr^n(r-1)}{2}.$$
One can see from (4) that all the coefficients of this polynomial are integers.
\qed

When $n=1$, the determinant in the left hand side of (2) is a common
Wronskian $W(1,f,\dots,f^{r-1})$. In this case $H_r$ is a polynomial
in one variable, which can be easily computed. The matrix
$W^{\prime}\in M(r,R)$ is a triangular matrix with the diagonal
elements $W_{kk}^{\prime}=k!(DF)^k,\,0\le k\le r-1$, hence $\det
W=\det W^{\prime}$ is equal to the product of these elements. So, we
have the following equality
$$\det\left\|D^kf^l\right\|_{0\le k,l\le r-1}=(Df)^{\frac{r(r-1)}{2}}
\prod_{k=1}^{r-1} k!\eqno{(5)}$$ where $f\in R$ and $D\in\Der(R)$.

The  formula (5) may be considered a consequence of the known
Wronskian chain rule [3,Part Seven, Ex. 56].

\section{ THE DETERMINANT}

For the rest of the paper $\ks$ is a field of fixed characteristic $p>0$, and
$q=\frac{p^n(p-1)}{2}$.
For $F\in\ks[X]^n$, denote
$$\Delta(F)=\det U(F)\in\ks[X^p].$$
Our aim is to compute this determinant.

Denote by $\phi_F$ the algebra endomorphism
$$\phi_F\in\End(\ks[X]),\,\phi_F:x_i\mapsto f_i,\,1\le i\le n.$$

{\bf Lemma 2}
{\it Let $F,G\in\ks[X]^n$. Let
$\phi_F G=(\phi_FG_1,\dots,\phi_FG_n)\in\ks[X]^n$.
Then
$$\Delta(\phi_FG)=(\phi_F\Delta(G))\cdot\Delta(F).\eqno{(6)}$$
}
\newline
{\it Proof}

By definition,
$$G^{\alpha}=\sum_{\gamma}U(G)_{\alpha\gamma}X^{\gamma},\,\alpha\in[0,p-1].$$
Applying the endomorphism $\phi_F$ to the both sides of this equality, we get
$$(\phi_FG)^{\alpha}=\sum_{\gamma}(\phi_FU(G)_{\alpha\gamma})F^{\gamma}=
\sum_{\beta\gamma}(\phi_FU(G)_{\alpha\gamma})U(F)_{\gamma\beta}X^{\beta}.$$
On the other hand,
$$(\phi_FG)^{\alpha}=\sum_{\beta}U(\phi_FG)_{\alpha\beta}X^{\beta}.$$
As the monomials $X^{\beta}$ form a basis, the coefficients are the same:
$$U(\phi_FG)=(\phi_FU(G))U(F).$$
Taking the determinant, we have (6). \qed

If $F\in\ks[X]^n$ consists of linear forms, then
$$F_i=\sum_{j=1}^n A_{ij}X_j,\,1\le i\le n$$
for some matrix $A\in M(n,\ks)$. We write this as
$$F=AX;$$
in this notation $X$ and $F$ are considered column vectors.

{\bf Lemma 3}
{\it Let $A\in M(n,\ks)$. If $F=AX$, then
$$\Delta(F)=(\det A)^{q}.\eqno{(7)}$$
}
\newline
{\it Proof}

Let $A,B\in M(n,\ks)$. Let $G=BX$ and $F=AX$. Then $\phi_FG=BAX$.
One can see that the elements of the matrix $U(BX)$ belong to the
field $\ks$. It follows that $\phi_F\Delta(BX)=\Delta(BX)$, and by
(6) we have
$$\Delta(BAX)=\Delta(BX)\Delta(AX).$$
It is well known that any matrix over a field is a product of diagonal matrices and
elementary ones. Thus, because of the latter equality,
it is sufficient to prove (7) for diagonal and elementary matrices.
If $A$ is elementary, then for some $j\neq k$ we have
$f_j=x_j+\lambda x_k$, where $\lambda\in\ks$, and
$f_i=x_i,\,i\neq j$. Then $U(F)$ is a (upper or lower) triangular matrix
with the unit diagonal, hence $\Delta(F)=\det A=1$, and (7) holds.
If $A$ is diagonal then $f_i=A_{ii}x_i,\,1\le i\le n$. In this case
$U(F)$ is diagonal as well and
$U(F)_{\alpha\alpha}=\prod_{i=1}^nA_{ii}^{\alpha_i}$. The equality (7)
can be easily checked.
\qed

{\bf Lemma 4}
{\it  Let $F\in\ks[X]^n$ and
$W=\left\|\partial^{\alpha}F^{\beta}\right\|_{\alpha,\beta\in[0,p-1]}$.
Then
$$\det W=
c_p^n\Delta(F),$$
where $c_p=\prod_{k=1}^{p-1}k!\in\ks^{\times}$.
}
\newline
{\it Proof}

Denote $Q=\left\|\partial^{\alpha}X^{\beta}\right\|
_{\alpha,\beta\in[0,p-1]}\in M(p^n,\ks[X])$. This is a Kronecker
product
$$Q=Q_1\otimes\dots\otimes Q_n;\,
Q_i=\left\|\partial_i^{a}x^{b}_i\right\| _{0\le a,b\le p-1}\in
M(p,\ks[x_i]),1\le i\le n.$$ Applying the equality (5) to the rings
$\ks[x_i]$, we have $\det Q_i=c_p, 1\le i\le n$. Then
$$\det Q=\prod_{i=1}^n(\det Q_i)^{p^{n-1}}=c_p^{np^{n-1}}=c_p^n.$$

The elements of $U(F)$ belong to the kernels of derivations $\partial_i$, hence
$$\partial^{\alpha}F^{\beta}=\sum_{\gamma}U(F)_{\beta\gamma}
\partial^{\alpha}X^{\gamma}.$$
This can be written in the matrix notation as
$$W=QU(F)^{\rm T}.$$ The formula is a consequence.
\qed

\section{PROOF OF THE THEOREM}

{\bf Proposition 2}
{\it
Let $F\in\ks[X]^n$. Then
$$\Delta(F)=j(F)^{q}.$$
}
\newline
{\it Proof}

Let $W=\left\|\partial^{\alpha}F^{\beta}\right\|
_{\alpha,\beta\in[0,p-1]}\in M(p^n,\ks[X])$.
By Lemma 4 and Proposition 1,
$$\det W=c_p^{n}\Delta(F)=H_{p}(JF).$$
If $A\in M(n,\ks)$ and $F=AX$ then $JF=A$, hence
$$H_{p}(A)=c_p^{n}(\det A)^{q}$$
by Lemma 3. Without loss of generality, $\ks$ is infinite.
The latter equality holds for any matrix, hence it is valid as a formal
equality in the ring $\ks[Z_{ij}]$. Thus
$$H_{p}(JF)=c_p^{n}(\det JF)^{q}=c_p^{n}j(F)^{q}.$$
\qed

We have the following corollary, proved first by Nousiainen [2].

{\bf Corollary}
{\it
Let $\ks$ be a field of characteristic $p>0$ and
$F\in\ks[X]^n$. Then the set
$\{F^{\alpha}:\alpha\in[0,p-1]\}$ is a basis of $\ks[X]$
over $\ks[X^p]$ if and only if $j(F)\in\ks^{\times}$.
}

{\bf Proposition 3}
{\it
Let $F\in\ks[X]^n$. Then
$$\ks[X]j(F)^{q}
\subset\ks[X^p][F].$$
}
\newline
{\it Proof}

From linear algebra we have
$$\Delta(F)U(F)^{-1}={\rm adj}\,U(F)\in M(p^n,\ks[X^p]).$$
Thus
$$j(F)^qX^{\alpha}=
\Delta(F)X^{\alpha}=\Delta(F)\sum_{\beta}U(F)^{-1}_{\alpha\beta}F^{\beta}
\in\ks[X^p][F]$$
for any multiindex $\alpha\in[0,p-1]$.
The set $\{X^{\alpha}\}$ is a basis of $\ks[X]$ over $\ks[X^p]$, hence
the inclusion follows.
\qed
\newline
{\it Proof of Theorem}

We have
$$J(R)=\langle\{j(F):F\in R^n\}\rangle=\langle P\rangle,\,P\in\ks[X].$$
If $P=0$, the statement is trivial. Suppose $P\neq 0$.
Since $P\in J(R)$, there exists a number $m\ge 1$, such that
$$P=\sum_{i=1}^mg_ij(F_i),$$
where $g_i\in\ks[X],\,F_i\in R^n,1\le i\le m$.
On the other hand, $J(R)\subset\langle P\rangle$, hence
$$\mu_i=j(F_i)/P\in\ks[X],\,1\le i\le m.$$

Consider the following two ideals:
$$J=\{f\in \ks[X]: \ks[X]f\subset R[X^p]\};\,
I=\langle \mu_1^q,\dots,\mu_m^q\rangle.$$
By Proposition 3, $\mu_i^qP^q=j(F_i)^q\in J$ for all $1\le i\le m$, hence
$$IP^q\subset J.$$
Raising the equality $\sum_{i=1}^mg_i\mu_i=1$ to the power $qm$, it
follows that $1\in I$, whence $P^q\in J$. \qed

\section*{ACKNOWLEDGMENT}

The author thanks Prof. Arno van den Essen for giving him the
information about the Nousiainen's preprint.

\section*{REFERENCES}

[1] H. Bass, E.H. Connell, D. Wright. The Jacobian Conjecture:
Reduction of degree and formal expansion of the inverse, Bull.
A.M.S. 7(2)(1982) 287-330.

[2] P. Nousiainen. On the Jacobian problem in positive
characteristic. Pennsylvania State Univ. preprint. 1981

[3] G. Polya, G. Szeg\"o. Problems and Theorems in Analysis, Vol II.
Springer-Verlag. 1976

\end{document}